\documentstyle[amscd,amssymb,verbatim,epsf]{amsart}

\begin{document}
\theoremstyle{plain}
\newtheorem{Ex}{Example}
\newtheorem{Thm}{Theorem}
\newtheorem{Cor}{Corollary}
\newtheorem{Main}{Main Theorem}
\newtheorem{Note}{Note}
\newtheorem{Lem}{Lemma}
\newtheorem{Prop}{Proposition}

\theoremstyle{definition}
\newtheorem{Def}{Definition}

\theoremstyle{remark}
\newtheorem{notation}{Notation}
\renewcommand{\thenotation}{}

\errorcontextlines=0
\numberwithin{equation}{section}
\renewcommand{\rm}{\normalshape}%

\title[Generalised Surfaces in ${\Bbb{R}}^3$]%
   {Generalised Surfaces in ${\Bbb{R}}^3$}
\author{Brendan Guilfoyle}
\address{Brendan Guilfoyle\\
          Department of Mathematics and Computing \\
          Institute of Technology, Tralee \\
          Clash \\
          Tralee  \\
          Co. Kerry \\
          Ireland.}
\email{brendan.guilfoyle@@ittralee.ie}
\thanks{The first author was supported by the Isabel Holgate Fellowship from  Grey College, Durham and the Royal Irish Academy Travel Grant Scheme}
\author{Wilhelm Klingenberg}
\address{Wilhelm Klingenberg\\
 Department of Mathematical Sciences\\
 University of Durham\\
 Durham DH1 3LE\\
 United Kingdom.}
\email{wilhelm.klingenberg@@durham.ac.uk }

\date{February 18, 2004}

\begin{abstract}
The correspondence between 2-parameter families of oriented lines in ${\Bbb{R}}^3$ and surfaces in $T{\Bbb{P}}^1$ is studied, and the 
geometric properties of the lines are related to the complex geometry of the surface. Congruences
generated by global sections of $T{\Bbb{P}}^1$ are investigated and a number
of theorems are proven that generalise results for closed convex
surfaces in ${\Bbb{R}}^3$.
\end{abstract}

\maketitle

\section{Introduction}

In this paper we study a generalisation of the concept of a
surface in Euclidean ${\Bbb{R}}^3$ through the twistor
construction. The fundamental aspect of this generalisation is that,
through its normal, an oriented surface in ${\Bbb{R}}^3$  determines a surface
in the space of all oriented lines: ${\Bbb{T}}=T{\Bbb{P}}^1$. This
4-manifold inherits a natural complex structure from the Euclidean
metric in ${\Bbb{R}}^3$ and has been useful in
connection with monopoles \cite{hitch1}.

A key feature of the construction is that 
a point in ${\Bbb{R}}^3$ corresponds to a holomorphic sphere in
${\Bbb{T}}$. In this paper we consider (not necessarily holomorphic)
2-sphere's in ${\Bbb{T}}$ and the {\it line congruences} they
generate in ${\Bbb{R}}^3$ \cite{pottwall}. Our main aim is to show that, in this
setting, curvature and umbilics can all be sensibly defined
and used to obtain results generalising those of closed convex surfaces.

The sense in which we are generalising is given by Frobenius' theorem:
we are allowing consideration of {\it twisting} line congruences. This
twist is encoded in the anti-symmetric part of a suitably defined second
fundamental form of the line congruence. The curvature $K$ of a congruence can then be defined as
the determinant of the second fundamental form.

This has the following bundle interpretation:

\noindent{\bf Theorem 1.}
{\it
A line congruence $\Sigma$ is locally the graph of a section of the bundle
$\pi:{\Bbb{T}}\rightarrow {\Bbb{P}}^1$ if and only if the curvature is
non-zero. 
}

We say a line congruence is {\it globally
convex} if it is the graph of a global section of the bundle
$\pi:{\Bbb{T}}\rightarrow {\Bbb{P}}^1$. Thus, a globally convex congruence $\Sigma\subset {\Bbb{T}}$ is a topological
sphere and the Gauss map $\pi|\Sigma$ yields natural global coordinates ($\xi$,$\bar{\xi}$) on $\Sigma$. We show that a generalised Gauss-Bonnet theorem
holds for globally convex congruences:

\noindent{\bf Theorem 2.}
{\it 
Let $\Sigma$ be a globally convex congruence with curvature
$K$. Then 
\[
\int_\Sigma Kd\mu=4\pi
\]
where $d\mu$ at $\gamma\in\Sigma$ is the pull-back via $\pi$ of the volume form induced on the
plane orthogonal to $\gamma$ by the Euclidean metric on ${\Bbb{R}}^3$.
}

In general these globally convex congruences, aside from twisting, will be
non-holomorphic. However, we can always perturb a congruence so that
it has isolated complex points. The total number of complex points,
counted with index, is a topological invariant of the congruence. We
show that this index is the number of shear-free points on a globally
convex congruence and that:

\noindent{\bf Theorem 3.}
{\it
The total number of shear-free lines (counted
with index) on a globally convex congruence with only isolated
shear-free lines is 4.
}

This generalises the well-known result that the number of isolated
umbilics (counted with index) on a closed convex surface is 4.

This paper is organised as follows: in the next section we recall the
relevant twistor construction, for further details see \cite{hitch1}
\cite{hitch2}.  To prove our main results we use a canonical
co-ordinate system, as in \cite{gak1} and null frame adapted to the congruence under
consideration. After describing these in section 3, we turn to the
first order description of congruences. Here the method of
spin-coefficients \cite{nap} \cite{perjes}, applied to ${\Bbb{R}}^3$,
yields a compact description of the geometric data. Finally in section
5 we prove the results regarding globally convex congruences.

\section{The Minitwistor Construction}

We begin by recalling the minitwistor construction of straight lines in ${\Bbb{R}}^3$ 
(see Hitchin \cite{hitch1} for further details).  Given a choice of origin in ${\Bbb{R}}^3$, 
a straight line can be uniquely described by two vectors: the (oriented)
direction of the line 
$\vec{W}$ and its perpendicular displacement from the origin $\vec{U}$. 
A straight line $\gamma$ is given by
\[
\gamma=\{\vec{U}+r\vec{W}\in {\Bbb{R}}^3\;|\; <\vec{U},\vec{W}>=0 \quad |\vec{W}|=1\quad r\in {\Bbb{R}}\},
\]
where $<,>$ is the Euclidean inner product and $|.|$ the associated norm.

The space of all oriented straight lines in ${\Bbb{R}}^3$ is the minitwistor space
\[
{\Bbb{T}}=\{(\vec{U},\vec{W})\in {\Bbb{R}}^3\times {\Bbb{R}}^3\;| \;<\vec{U},\vec{W}>=0 \quad |\vec{W}|=1 \}
\cong TS^2.
\]
Throughout this paper we utilise this bijection to identify an oriented line $\gamma\subset{\Bbb{R}}^3$ with the point $\gamma\in{\Bbb{T}}$. This 4-manifold has a natural almost complex structure ${\Bbb{J}}$ defined by
rotation in ${\Bbb{R}}^3$
through 90$^o$ about the direction of the line.  In fact, the almost complex 
structure ${\Bbb{J}}$ is integrable and so $\Bbb{T}$ is a complex surface.  Alternatively, ${\Bbb{J}}=j\oplus j$, where the splitting $T{\Bbb{T}}=TS^2\oplus TS^2$ and the complex structure $j$ on $S^2$ are induced by the Euclidean metric on ${\Bbb{R}}^3$.  In addition, 
there is an anti-holomorphic involution $\tau: \Bbb{T} \rightarrow \Bbb{T}$
given by reversing the orientation of the line.

A point $p$ in ${\Bbb{R}}^3$ is uniquely determined by the 2-sphere
$S^2_p$ of oriented  lines 
passing through it.  This sphere $S^2_p\subset \Bbb{T}$ has the following
properties:
\begin{enumerate}
 \item $S^2_p$ is a complex line in $ \Bbb{T}$ i.e. the complex
 structure on $ \Bbb{T}$  leaves invariant the tangent space of $S^2_p$
 \item $S^2_p\cap S^2_q$ consists of two points in $\Bbb{T}$ - the two oriented lines in ${\Bbb{R}}^3$ passing through $p$ and $q$
 \item $S^2_p$ is invariant under the involution $\tau$
\end{enumerate}

From the above it follows that $S^2_p$ can be given, in terms of a
holomorphic coordinate $\xi$ on ${\Bbb{P}}^1$, as the graph of a section of $\pi:{\Bbb{T}}
\rightarrow {\Bbb{P}}^1$
\begin{equation}\label{e:cp1}
s(\xi)= \frac{1}{2}\left((x^1+ix^2)-2x^3\xi-(x^1-ix^2)\xi^2\right)\frac{\partial}{\partial \xi},
\end{equation}
where the Euclidean coordinates of $p$ are ($x^1,x^2,x^3$).  This
follows from the fact that global holomorphic sections can be at most
quadratic in $\xi$ (self-intersection 2) and then the invariance of
the section under the antipodal map $\tau(\xi)= -\overline{\xi}^{-1}$
restricts the coefficients to be of the above form.

In this paper we will investigate {\it line congruences}, that is,
surfaces $\Sigma\subset \Bbb{T}$ or two parameter families of lines in
${\Bbb{R}}^3$. Every surface $S\subset{\Bbb{R}}^3$ gives rise to a
congruence $\Sigma\subset \Bbb{T}$ by way of its normal, but not every
congruence arises in this way. It is in this sense that we are working
with generalised surfaces.

In the next section we introduce new local coordinates and a null frame on
${\Bbb{R}}^3$ which fit nicely with the description of a congruence as a
surface in $T{\Bbb{P}}^1$.

\section{Coordinates and an Adapted Frame on ${\Bbb{R}}^3$}
Let $(x^1,x^2,x^3)$ be the standard
coordinates on ${\Bbb{R}}^3$ and and set $z=x^1+ix^2$, $\overline{z}=x^1-ix^2$ and $x^3=t$.

\begin{Def} Consider
the following transformation $ \Phi:(u,v,r)\rightarrow
(z,\overline{z},t)$ on an open subset of ${\Bbb{R}}^3$ given by 
\begin{equation}\label{e:coord1}
z=\frac{2(F-\overline{F}\xi^2)+2\xi(1+\xi\overline{\xi})r}{(1+\xi\overline{\xi})^2}
\end{equation}
\begin{equation}\label{e:coord3}
t=\frac{-2(F\overline{\xi}+\overline{F}\xi)+(1-\xi^2\overline{\xi}^2)r}{(1+\xi\overline{\xi})^2},
\end{equation}
where $F(u,v)$ and $\xi(u,v)$ are smooth complex-valued functions of
two real parameters $u$ and $v$.  We call the coordinates ($u,v,r$) the {\it congruence  coordinates}. 
\end{Def}

For each ($u,v$), $\Phi(r)$ is a straight line in ${\Bbb{R}}^3$.
Moreover, these lines are parameterised by arclength $r$ and
$\frac{\partial}{\partial r}$ is the unit tangent to the lines. The
shortest distance from the origin to each line is given by the point
$r=0$. 

These coordinates come from the twistor construction in the following way:  consider one of these
lines with direction $\vec{W}$ and perpendicular displacement from the origin $\vec{U}$.  Translate
$\vec{W}$ along $\vec{U}$ to the origin and then translate $\vec{U}$
along $\vec{W}$.  This vector tangent to $S^2$ gives us the
line as a point in $\Bbb{T}$.  Moreover, $\xi$ is the standard holomorphic coordinate on ${\Bbb{P}}^1$ induced from 
stereographic projection from the South pole and in these coordinates the point in $\Bbb{T}$ is given by
\[
\left(\xi(u,v),F(u,v)\frac{\partial}{\partial
    \xi}\right)\in T_{\xi}{\Bbb{P}}^1.
\]
Thus $F(u,v)$ determines the perpendicular distance of the line from
the origin. This can be viewed as parametric equations for a line congruence in terms of coordinates ($\xi$,$F$) on ${\Bbb{T}}$, which are holomorphic with respect to ${\Bbb{J}}$. For further details see \cite{gak1}.

A change of origin leads to a quadratic holomorphic
translation of the function $F$.  In particular, if the origin is translated 
$(0,0,0) \rightarrow (x^1_0,x^2_0,x^3_0)$ then 
\begin{equation}\label{e:trans1}
F\rightarrow F+\frac{1}{2}\left(\alpha_0-2t_0\xi-\overline{\alpha}_0\xi^2\right),
\end{equation}
where $\alpha_0=x^1_0+ix^2_0$ and $t_0=x_0^3$. This can be seen from
equation (\ref{e:cp1}) since the lines through the origin (the zero
section of the bundle) will change to the section
\[
s(\xi)= \frac{1}{2}\left(\alpha_0-2t_0\xi-\overline{\alpha_0}\xi^2\right)\frac{\partial}{\partial \xi}.
\]
 Our coordinates will then transform by:
\begin{equation}\label{e:trans3}
(u,v,r)\rightarrow \left(u,v,r+\frac{\overline{\alpha}_0\xi+\alpha_0\overline{\xi}+t_0(1-\xi\overline{\xi})}
                         {1+\xi\overline{\xi}}\right)
\end{equation}
\[
(z,t)\rightarrow(z+\alpha_0,t+t_0).
\]
The change in $r$ is just $<\vec{T},\vec{W}>$, where $\vec{T}$ is
the translation vector determined by $\alpha_0$ and $t_0$.

The quantities that have geometric significance are invariant under
this translation.  In particular, we have the following translation
invariant derivatives of the twistor function $F$:

\begin{Prop}\label{p:traninv}
Let $\nu=u+iv$,$\overline{\nu}=u-iv$, and $\partial =\frac{\partial
}{\partial\nu}$ and $\overline{\partial}=\frac{\partial
}{\partial\overline{\nu}}$. Then 
\[
\partial^+F\equiv\partial F+r\partial\xi-\frac{2F\overline{\xi}\partial \xi}{1+\xi\overline{\xi}}
\]
\[
\partial^-F\equiv\overline{\partial} F+r\overline{\partial}\xi-\frac{2F\overline{\xi}\;\overline{\partial} \xi}{1+\xi\overline{\xi}},
\]
are invariant under the translations (\ref{e:trans1}) and (\ref{e:trans3}).

\end{Prop}
\begin{pf}
This is a straight-forward calculation. 

\end{pf}

\begin{Note}
The Jacobian of the tranformation $\Phi$ is
\[
\Delta=\frac{4}{(1+\xi\overline{\xi})^2}\left(\partial^+F\overline{\partial^+F}-\partial^-F\overline{\partial^-F}\right)
\]
Thus the transformation is a diffeomorphism wherever $\Delta\neq0$.
\end{Note}

A {\it null frame} in ${\Bbb{R}}^3$ is a trio $\{e_0,e_+,e_-\}$ of
complex vector fields in ${\Bbb{C}}\otimes T\:{\Bbb{R}}^3$, where $e_0$ is real, 
$e_+$ is the complex conjugate of $e_-$ and they satisfy the following orthogonality properties:
\[
<e_0, e_0>=1 \qquad <e_0, e_+>=0 \qquad <e_+, e_+>=0 \qquad <e_+,
e_->=1,
\]
where we have extended the Euclidean inner product of ${\Bbb{R}}^3$ bilinearly over
${\Bbb{C}}$. Orthonormal frames $\{e_0,e_1,e_2\}$ on $T\:{\Bbb{R}}^3$ and null frames are
related by
\begin{equation}\label{e:ortho}
e_+=\frac{1}{\sqrt{2}}(e_1-ie_2) \qquad
e_-=\frac{1}{\sqrt{2}}(e_1+ie_2).
\end{equation}

\begin{Def}
A {\it congruence null frame} for $\Sigma\subset{\Bbb{T}}$ is a null frame
$\{e_0,e_+,e_-\}$ if, for each $\gamma\in\Sigma$, we have $e_0$ tangent
to $\gamma$ in ${\Bbb{R}}^3$, and the orientation of $\{e_0,e_1,e_2\}$
is the standard orientation on ${\Bbb{R}}^3$.
\end{Def}

\begin{Prop}\label{p:nullframe}
Let $\Sigma$ be a line congruence and consider an open set $U\subset\Sigma$ with $\Delta\neq 0$. Suppose that $F$ is the twistor function describing
$\Sigma$ on $U$.
A null frame is a congruence null frame to $\Sigma$ if and only if it
has the following expression in terms of canonical coordinates 
\[
e_0=\frac{\partial}{\partial r}
\qquad\qquad
e_+= \left(\alpha\frac{\partial}{\partial \nu}
     +\beta\frac{\partial}{\partial \overline{\nu}}
     +\Omega\frac{\partial}{\partial r}\right)e^{i\phi}
\qquad\qquad
e_-=\overline{e_+},
\]
where
\[
\Omega=\frac{\sqrt{2}\left[\overline{\partial^-F}(F\overline{\partial}\;\overline{\xi}+\overline{F}\;\overline{\partial}\xi)-\overline{\partial^+F}(F\partial\overline{\xi}+\overline{F}\partial\xi)\right]}
{(1+\overline{\xi}\xi) (\partial^-F\overline{\partial^-F}-\partial^+F\overline{\partial^+F}) }
\]
\begin{equation}\label{e:coef}
\alpha=\frac{\overline{\partial^+F}(1+\overline{\xi}\xi)}
{\sqrt{2}(\partial^-F\overline{\partial^-F}-\partial^+F\overline{\partial^+F}) }
\qquad
\beta=-\frac{\overline{\partial^-F}(1+\overline{\xi}\xi)}
{\sqrt{2}(\partial^-F\overline{\partial^-F}-\partial^+F\overline{\partial^+F}) },
\end{equation}
and $\phi$ is a function of $\nu$, $\overline{\nu}$ and $r$.
\end{Prop}
\begin{pf}
We can prove this by a change of coordinates from
 $(\nu,\overline{\nu},r)$ to $(z,\overline{z},t)$ using equations
 (\ref{e:coord1}) to (\ref{e:coord3}). This gives the following mixed
 form for the frame:
\[
e_0= \frac{2\xi}{1+\xi\overline{\xi}}\frac{\partial}{\partial z}
     +\frac{2\overline{\xi}}{1+\xi\overline{\xi}}\frac{\partial}{\partial \overline{z}}
  + \frac{1-\xi\overline{\xi}}{1+\xi\overline{\xi}}\frac{\partial}{\partial t}
\]
\[
e_+= \frac{\sqrt{2}}{1+\xi\overline{\xi}}\frac{\partial}{\partial z}
     -\frac{\sqrt{2}\;\overline{\xi}^2}{1+\xi\overline{\xi}}\frac{\partial}{\partial \overline{z}}
  - \frac{\sqrt{2}\;\overline{\xi}}{1+\xi\overline{\xi}}\frac{\partial}{\partial t}.
\]
Now, since the Euclidean inner product on ${\Bbb{R}}^3$ in the coordinates
$z,\overline{z},t$ is simply 
\[
g_{ij}=\left[\begin{matrix}
                 0 & \frac{1}{2} & 0\\
                 \frac{1}{2} & 0 & 0\\
                 0 & 0 & 1 
              \end{matrix}\right],
\]
we can check that these
vectors form a null frame. 

Alternatively, we can work backwards as follows. Set 
\[
e_0=\frac{\partial}{\partial r}
\qquad\qquad
e_+= \alpha\frac{\partial}{\partial \nu}
     +\beta\frac{\partial}{\partial \overline{\nu}}
     +\Omega\frac{\partial}{\partial r},
\]
for some functions $\alpha$, $\beta$ and $\Omega$ to be determined. Now the first condition on these functions is $<e_0, e_+>=0$ which gives us that
\[
\Omega =-\alpha\;\left<e_0,\frac{\partial}{\partial\nu}\right>- \beta\left<e_0,\frac{\partial}{\partial\overline{\nu}}\right>.
\]
Thus
\begin{equation}
e_+=\alpha\left(\frac{\partial}{\partial\nu}-\left<e_0,\frac{\partial}{\partial\nu}\right>\frac{\partial}{\partial r}\right)
 +\beta\left(\frac{\partial}{\partial\overline{\nu}}-\left<e_0,\frac{\partial}{\partial\overline{\nu}}\right>\frac{\partial}{\partial r}\right)
\quad=\alpha Z_++\beta Z_-\nonumber.
\end{equation}
The second condition $<e_+, e_+>=0$ now becomes
\begin{equation}\label{e:quad1}
<Z_+, Z_+>\alpha^2+2<Z_+, Z_->\alpha\beta\;+<Z_-, Z_->\beta^2=0.
\end{equation}
A lenghty calculation, involving a change to Euclidean coordinates via
$\Phi$, shows that
\begin{equation}\label{e:inprod}
 <Z_+, Z_+>=\frac{4\partial^+F\overline{\partial^-F}}{(1+\xi\overline{\xi})}
\qquad
<Z_+,Z_->=\frac{2(\partial^-F\overline{\partial^-F}+\partial^+F\overline{\partial^+F})}
 {(1+\xi\overline{\xi})}.
\end{equation}
With the aid of this, we can solve equation (\ref{e:quad1}) for the 
ratio of $\alpha$ and $\beta$:
\[
\frac{\alpha}{\beta}=-\frac{\overline{\partial^+F}}{\overline{\partial^-F}}
 \qquad or \qquad 
\frac{\alpha}{\beta}=-\frac{\partial^-F}{\partial^+F}.
\]
The fact that we get two answers simply represents the choice of
orientation for the null frame. We choose the first in order to have
our orientation agree with the standard orientation of
${\Bbb{R}}^3$. This orientation coincides with the graph orientation
$\Sigma$ in ${\Bbb{T}}\rightarrow{\Bbb{P}}^1$.

The final equation we have to solve is $<e_+, e_->=1$ and this gives
us $\beta\overline{\beta}$, that is $\phi$.  This represents rotation
of the frame  about $e_0$, which we can set to agree with the argument
of $-\overline{\partial^-F}$ by parallel translation of the frame.  
 
\end{pf}

\section{The Local Geometry of Congruences}

We now use the twistor coordinates and null frame of the last section
to investigate the local geometry of congruences in ${\Bbb{R}}^3$.
To start we find expressions for the first order geometric properties:

\begin{Prop}\label{p:opscal}
The Lie derivative of $e_+$ in the $e_0$ direction is given by
\[
L_{e_0}e_+=\overline{\rho}e_++\sigma e_-,
\]
with
\begin{equation}\label{e:rho}
\rho=\frac{ \partial^+F\overline{\partial}\;\overline{\xi} -\partial^-F\partial\overline{\xi}}
{\partial^-F\overline{\partial^-F}-\partial^+F\overline{\partial^+F}}
\end{equation}
\begin{equation}\label{e:shear}
\sigma=\frac{\overline{\partial^+F}\partial\overline{\xi} -\overline{\partial^-F}\;\overline{\partial}\;\overline{\xi}}
{\partial^-F\overline{\partial^-F}-\partial^+F\overline{\partial^+F}}.
\end{equation}
\end{Prop}

\begin{pf}
This follows by using the computations in the previous Proposition.  For the twist and divergence,
\begin{align}\nonumber
\overline{\rho} &= \left<L_{e_0}e_+, e_-\right> \\\nonumber
  &= \left<\frac{\partial \alpha}{\partial r}Z_++\frac{\partial \beta}{\partial r}Z_-\;,\;\overline{\alpha}Z_-+\overline{\beta}Z_+\right>\nonumber.
\end{align}
Now use the expressions for $\alpha$ and $\beta$ and equation (\ref{e:inprod}) to get the result.

Similarly
\begin{align}\nonumber
\sigma &= \left<L_{e_0}e_+, e_+\right> \\\nonumber
  &= \left<\frac{\partial \alpha}{\partial r}Z_++\frac{\partial \beta}{\partial r}Z_-\;,\;\alpha Z_++\beta Z_-\right>\nonumber.
\end{align}

\end{pf}

The complex scalar functions $\rho$ and $\sigma$ describe the first order geometric behaviour of the congruence 
of lines. In particular, the real part of $\rho$ is the {\it divergence}, the imaginary part is the 
{\it twist} and $\sigma$ is the {\it shear} of the congruence (see
\cite{guil} for details). By Proposition \ref{p:traninv} these are invariant under translations of
the origin.  For line congruences in ${\Bbb{R}}^3$ the
evolution of these quantities along the line can be more directly
derived.  In particular, they satisfy the {\it Sachs equations} \cite{par}
\begin{equation}\label{e:sachs}
\frac{\partial \rho}{\partial r}=\rho^2+\sigma\overline{\sigma}
\qquad 
\frac{\partial \sigma}{\partial r}=(\rho+\overline{\rho})\sigma.
\end{equation}
Note that if the shear or the twist vanish at some point on the line, they
vanish at every point on the line. 

Following Penrose and Rindler \cite{par} these can be integrated in
terms of the initial values of $\rho$ and $\sigma$ at $r=0$:

\[
\rho=\frac{\rho_0-(\rho_0\overline{\rho}_0-\sigma_0\overline{\sigma}_0)r}
  {1-(\rho_0+\overline{\rho}_0)r+(\rho_0\overline{\rho}_0
                              -\sigma_0\overline{\sigma}_0)r^2}
\qquad
\sigma=\frac{\sigma_0}
  {1-(\rho_0+\overline{\rho}_0)r+(\rho_0\overline{\rho}_0-\sigma_0\overline{\sigma}_0)r^2}.
\]

The shear has the following interpretation. Consider a circle in the plane orthogonal to the line $\gamma_0$ at some
point. Lie-propogation along the line can alter this circle in a number of
ways: if the shear is zero, it will remain a circle and
if the shear is non-zero it will become an ellipse. In particular,
$|\sigma|$ measures the eccentricity of the ellipse, while
$\phi=Arg(\sigma)$ measures the inclination of the semi-major and
semi-minor axes (details can be found in \cite{guil}).

The twist can be understood as follows:
\begin{Def} 
A line congruence is {\it integrable} iff locally there exists an embedded surface $S$ in ${\Bbb{R}}^3$ such that $S$ is orthogonal to the lines of the congruence.
\end{Def}

\begin{Prop}\cite{guil}
A line congruence is integrable iff $\rho$ is real (the twist vanishes).
\end{Prop}

\begin{Prop}
The orthogonal surface $S$ to an integrable line congruence is given by $r=r(\nu,\bar{\nu})$, where
\[
\bar{\partial} r=\frac{2F\bar{\partial}\bar{\xi}+2\bar{F}\bar{\partial}\xi}{(1+\xi\bar{\xi})^2}.
\]
\end{Prop}

\begin{pf}
It is not hard to show that the parametric surface in ${\Bbb{R}}^3$,
obtained by inserting $r=r(\nu,\bar{\nu})$ in (\ref{e:coord1}) and (\ref{e:coord3}),  is orthogonal to the
congruence iff the above condition holds. 
\end{pf}

\begin{Def}
The {\it curvature} of a congruence is defined to be $K=\rho\overline{\rho}-\sigma\overline{\sigma}$.
\end{Def}

This definition reflects the fact that, in the
twist-free case, the curvature of a congruence is the curvature of the one
parameter family of  surfaces in ${\Bbb{R}}^3$
orthogonal to the congruence. More generally,

\begin{Thm}\label{t:section}
A line congruence is locally the graph of a section of the bundle
$\pi:{\Bbb{T}}\rightarrow {\Bbb{P}}^1$ if and only if the curvature is
non-zero. 
\end{Thm}
\begin{pf}
From equations (\ref{e:rho}) and (\ref{e:shear}) we find that
\[
K=\frac{\partial\xi\overline{\partial}\;\overline{\xi}-\partial\overline{\xi}\;\overline{\partial}\xi}
{\partial^-F\overline{\partial^-F}-\partial^+F\overline{\partial^+F}}.
\]
Thus, the left hand side will vanish at a point $\gamma\in\Sigma$ if
and only if the Jacobian relating $\xi$ and $\nu$ is zero at $\gamma$,
that is, the tangent space to $\Sigma$ at $\gamma$ has a
vertical component. 
\end{pf}

\section{Globally Convex Congruences}

We say a congruence is {\it globally convex} if it is the graph of a global
section of the canonical bundle $\pi:{\Bbb{T}}\rightarrow
{\Bbb{P}}^1$. A globally convex congruence is a toplogical 2-sphere which generalises the concept of an closed convex surface. 

\vspace{0.1in}
\noindent{\bf Example}: The derivative of the action of $PSL(2,{\Bbb{C}})$ on ${\Bbb{P}}^1$
generates a holomorphic vector field on ${\Bbb{P}}^1$ which is a
global section of the canonical bumdle. This 6-parameter family of line congruences splits into twisting and twist-free congruences. The former 
contains the standard overtwisted contact structure on ${\Bbb{R}}^3$, see
\cite{elia}, while the latter consists of the holomorphic ${\Bbb{P}}^1$
generated by lines through a point in ${\Bbb{R}}^3$.

\vspace{0.1in}

Globally convex congruences have nice properties. For example,

\begin{Prop}
Suppose that $\Sigma\subset{\Bbb{T}}$ is a globally convex congruence,
then every
point in ${\Bbb{R}}^3$, is contained in some line of the congruence,
i.e. the congruence is spacefilling.
\end{Prop}
\begin{pf}
Consider any point $p\in {\Bbb{R}}^3$.  By a
translation we can make $p$ the origin.  Now the resulting surface
$\Sigma$ is a vector field on $S^2$ and therefore must have a zero.  In terms of the lines on ${\Bbb{R}}^3$, a  
zero represents a line with zero perpendicular distance from the origin, that is a line passing through the origin.
Thus every point $p$ must have a line passing through it.

\end{pf}

\begin{Thm}
Let $\Sigma$ be a globally convex congruence with curvature
$K$. Then:
\[
\int_\Sigma Kd\mu=4\pi,
\]
where $d\mu$ at $\gamma\in\Sigma$ is the pull-back via $\pi$ of the volume form induced on the
plane orthogonal to $\gamma$ by the Euclidean metric on ${\Bbb{R}}^3$.
\end{Thm}

\begin{pf}
Suppose we choose coordinates $(\xi,\overline{\xi})$ and twistor
function $F$ to describe an open subset of $\Sigma$.
 
The basis of 1-forms $\{\theta^0,\theta^+,\theta^-\}$ dual to the vector
basis in Proposition \ref{p:nullframe} has coordinates
\[
\theta^+=\frac{\sqrt{2}}{K(1+\xi\overline{\xi})}\left(\rho
  d\xi-\overline{\sigma} d\overline{\xi}\right).
\]
Now, $d\mu=i\pi^*(\theta^+\wedge\theta^-)$ and so we get that
\begin{equation}
\int_\Sigma K d\mu=\int_\Sigma K\;i\pi^*(\theta^+\wedge\theta^-)\quad
              =2\int_{{\Bbb{P}}^1}\frac{d\xi
              d\overline{\xi}}{(1+\xi\overline{\xi})^2}\quad
              =4\pi \nonumber.
\end{equation}

\end{pf}

Consider a line congruence in ${\Bbb{R}}^3$ given by an oriented surface
$\Sigma\subset{\Bbb{T}}$.  

\begin{Def} A point $\gamma\in\Sigma$ is {\it complex} if the complex
structure acting on ${\Bbb{T}}$ preserves $T_\gamma\Sigma$, the
tangent space to $\Sigma$ at $\gamma$. A complex point $\gamma\in\Sigma$ is
called {\it positive} or {\it negative} if the complex structure
induced on $T_\gamma\Sigma$ by ${\Bbb{J}}$ agrees or disagrees with
the given orientation on $\Sigma$, respectively. 
\end{Def}

In the general case a line congruence $\Sigma\subset{\Bbb{T}}$ will
have places where it is complex and others where it fails to
be. We can perturb such a surface so that the complex
points are isolated. Moreover, there is a well-defined index for each 
isolated complex point and we identify these as follows:
 
\begin{Prop}\label{p:minikerr}
A point on $\Sigma$ is complex iff the shear vanishes along the
line. The index $I$ of an isolated complex point $\gamma_0$ is minus the
winding number of the semi-major axes of shear as we go around $\gamma_0$.
\end{Prop}

\begin{pf}
Let $f:\Sigma\rightarrow {\Bbb{T}}$ be a $C^2$ smooth immersion with
differential $df:T\Sigma\rightarrow T{\Bbb{T}}$.  Let $j$ be a
conformal structure on $\Sigma$ compatible with $\nu$. We define the
sections $\delta^+f\in\Omega^{10}(\Sigma)\otimes T^{10}{\Bbb{T}}$ and 
$\delta^-f\in\Omega^{01}(\Sigma)\otimes T^{10}{\Bbb{T}}$ by
\[
\delta^\pm f=\frac{1}{2}\left(df \mp {\Bbb{J}}\circ df \circ j\right).
\]
Then
$\delta^+f\wedge\delta^-f\in\Omega^2(\Sigma)\otimes\mbox{det}\;T^{10}{\Bbb{T}}$,
works out to be 
\[
\delta^+f\wedge\delta^-f=(\partial\xi\bar{\partial}\eta-\bar{\partial}\xi\partial\eta)d\nu\wedge d\bar{\nu}
\]

A point is complex iff this 2-form vanishes, see \cite{kling}. By
simplifying the numerator of $\sigma$ in (\ref{e:shear}), this is
equivalent to
the vanishing of the shear.

The index of a complex point is equal to the intersection index of
$\delta^+f\wedge\delta^-f$ with the zero section of
$\Omega^2(\Sigma)\otimes\mbox{det}\;T^{10}{\Bbb{T}}$. This is just the
winding number of the complex function
$\partial\xi\bar{\partial}\eta-\bar{\partial}\xi\partial\eta$,
i.e. the winding number of $\bar{\sigma}$ about its isolated zero. 
\end{pf}

In the integrable case the argument of the shear is the angle between
the real axis of $\nu$ with a principal curvature direction of the
orthogonal surface in
${\Bbb{R}}^3$. We have the following generalisation of the fact that the total number
of isolated umbilics (counted with index) on an closed convex surface is 4.

\begin{Thm}
The total number of shear-free lines (counted
with index) on a globally convex congruence with only isolated
shear-free lines is 4.
\end{Thm}
\begin{pf}

Let $d_+$ be the sum of the indices over all positive complex points
of $\Sigma$. Define $d_-$ similarily for negative complex points. Then
by \cite{lai},
\[
d_++d_-=\chi(T\Sigma)+\chi(N\Sigma)
\]
\[
d_+-d_-=c_1(TM)\left[\Sigma\right],
\]
where $N\Sigma$ is the normal bundle and $T\Sigma$ is
the tangent bundle of $\Sigma$,
$\chi$ is the Euler number of the appropriate bundle, and $c_1(TM)$ is
the first Chern class of the tangent bundle to $M$.

In the case of a globally convex congruence with the graph orientation, we
have $\chi(T\Sigma)=2$, $\chi(N\Sigma)=2$ and
$c_1(T{\Bbb{T}})\left[\Sigma\right]=4$. The theorem then follows.

\end{pf}

\end{document}